\documentclass[11pt]{amsart}
\usepackage[T1]{fontenc}
\usepackage{lmodern}
\usepackage{amsmath,amssymb,mathrsfs,mathtools}
\usepackage[margin=1.15in]{geometry}
\usepackage{microtype}
\usepackage{needspace}
\usepackage{cite}
\usepackage[colorlinks=true,linkcolor=blue,citecolor=blue,urlcolor=blue]{hyperref}
\hypersetup{pdftitle={The Cp conjecture of Muckenhoupt -- Deterministic proof}}
\newtheorem{theorem}{Theorem}[section]
\newtheorem{proposition}[theorem]{Proposition}
\newtheorem{lemma}[theorem]{Lemma}
\newtheorem{corollary}[theorem]{Corollary}
\theoremstyle{remark}
\newtheorem{remark}{Remark}[section]
\numberwithin{equation}{section}
\newcommand{\R}{\mathbb R}
\newcommand{\D}{\mathscr D}
\newcommand{\Sfam}{\mathscr S}
\newcommand{\Pfamily}{\mathscr P}
\newcommand{\Ofamily}{\mathscr O}
\newcommand{\one}{\mathbf 1}
\newcommand{\avg}[2]{\frac{1}{|#2|}\int_{#2}#1}
\DeclareMathOperator*{\essinf}{ess\,inf}
\title{The $C_p$ conjecture of Muckenhoupt}

\author{Andrei K. Lerner}
\address{Department of Mathematics,
Bar-Ilan University, 5290002 Ramat Gan, Israel}
\email{lernera@math.biu.ac.il}

\author{Sheldy Ombrosi}
\address{Departamento de Matem\'atica\\
Universidad Nacional del Sur\\
Bah\'ia Blanca, 8000, Argentina}\email{sombrosi@uns.edu.ar}

\begin{document}
\begin{abstract}
We give an affirmative answer to Muckenhoupt's $C_p$ conjecture:
for $1<p<\infty$, the $C_p$ condition is necessary and sufficient
for the weighted $L^p(w)$ norm inequality relating the Hilbert
transform and the Hardy--Littlewood maximal operator. We also
show that the $C_p$ condition is necessary and sufficient for
the weighted Fefferman--Stein inequality.
\end{abstract}


\maketitle

\section{Introduction}

The main results proved in this paper hold in every dimension.
For clarity of exposition, we present them
in the one-dimensional setting; the extension to $\R^n$ is
explained in Remark~\ref{rem:higher-dimensions}.

Let $H$ and $M$ denote the Hilbert transform and the uncentered
Hardy--Littlewood maximal operator, respectively. For $1<p<\infty$,
a weight $w$ on $\R$ belongs to $C_p$ if there exist $r>1$ and
$K<\infty$ such that
\[
 \left(\avg{w^r}{I}\right)^{1/r}
 \leq \frac{K}{|I|}\int_{\R}(M\one_I)^p w
\qquad\text{for every bounded interval }I.
\]

Muckenhoupt~\cite{Muckenhoupt} proved that the $C_p$ condition is necessary
for the inequality
\begin{equation}\label{eq:main}
 \int_{\R}|Hf|^p w\leq C\int_{\R}(Mf)^p w
\end{equation}
and conjectured that it is also sufficient. This conjecture has
remained open for 46 years.

Sawyer~\cite{Sawyer} established sufficiency under the
stronger assumption $w\in C_{p+\varepsilon}$ for some $\varepsilon>0$.
This raised the question whether the $C_p$ condition is self-improving
in $p$: an affirmative answer would have settled the conjecture.
However, Kahanp\"a\"a and Mejlbro~\cite{KM} constructed a
weight in $C_p$ which belongs to no $C_q$ with $q>p$, ruling out
this approach.

Further developments on $C_p$ weights include reverse H\"older estimates, extensions to rough and vector-valued operators, and characterizations of weak-type inequalities; see
\cite{Buckley,Canto,CantoLiRoncalTapiola,CejasLiPerezRivera,LernerFS,LernerWeak,Yabuta}. Nevertheless, the original sufficiency problem remained unresolved.

Sparse domination, developed during the 2010s, allows us to replace
the Hilbert transform by simpler positive operators built from local
averages over sparse families of intervals. However, for more than
a decade it remained unclear how to prove \eqref{eq:main} with
these sparse operators in place of $H$, assuming only $w\in C_p$.
The difficulty lies in controlling the overlapping nonlocal tails
produced by the $C_p$ condition, even when the underlying family is
sparse. We establish this bound and thereby prove the Muckenhoupt
conjecture.

A family $\Sfam$ of dyadic intervals is
$\eta$-sparse, $0<\eta<1$, if there are pairwise disjoint measurable sets
$E_I\subset I$, $I\in\Sfam$, with $|E_I|\geq\eta|I|$. Its sparse operator is
\[
 A_{\Sfam}f=\sum_{I\in\Sfam}f_I\one_I,
 \qquad f_I=\avg{f}{I}.
\]

\begin{theorem}\label{thm:main}
Let $1<p<\infty$, let $w$ satisfy the $C_p$ condition with parameters
$r,K$, and let $\Sfam$ be an $\eta$-sparse family in a dyadic grid.
Then, for every locally integrable function $f$,
\[
 \int_{\R}|A_{\Sfam}f|^p w
 \leq C_{p,r,K,\eta}\int_{\R}(Mf)^p w.
\]
\end{theorem}

By the standard sparse domination of the Hilbert transform,
see~\cite{LN}, Theorem~\ref{thm:main} immediately gives the following.

\begin{corollary}\label{cor:hilbert}
Muckenhoupt's conjecture is true. More precisely, if $1<p<\infty$ and
$w\in C_p$, then \eqref{eq:main} is valid for every bounded compactly
supported function $f$, with a constant depending only on $p,r,K$.
\end{corollary}

\subsection{Idea of the proof}

We outline the main ideas of the proof, leaving the technical
details to the later sections. Write $w(I)=\int_Iw$ and $w_I=w(I)/|I|$,
and assume that $f\geq0$.

\medskip
\noindent\emph{From the sparse operator to a sum.}
For any weight $w$ and any $r>1$, Lemma~\ref{lem:sparse} gives
\[
 \int_{\R}(A_{\Sfam}f)^pw
 \lesssim
 \sum_{I\in\Sfam}f_I^p
       \left(\avg{w^r}{I}\right)^{1/r}|I|.
\]
The $C_p$ condition is not used at this stage. This lemma is the
main analytic use of the sparseness of the original family.

The natural next step would be to bound the sum on the right by
$\int_{\R}(Mf)^pw$. Under the $C_p$ condition alone, however, this
is false. Appendix~\ref{app:sparse-sum} gives a weight in $C_p$ and
a bounded compactly supported $f$ for which the maximal-function
integral is finite but this sum is infinite for every $r>1$.
We must therefore split the operator before using the lemma.
The intervals along which the maximal averages of $f$ grow rapidly
are handled directly by a pointwise estimate; all the intervals in
the counterexample belong to this part. We apply the lemma to the
remaining family.

For this remaining sum, the $C_p$ condition bounds the $L^r$ average
of $w$ by the nonlocal quantity
$|I|^{-1}\int_{\R}(M\one_I)^pw$. Direct summation of these
overlapping tails does not give the desired estimate. Our plan is
to transform the sum by telescoping. For this we need two
structures: monotonicity in the quantities involving $f$, and an
exact relation between the weight terms at a parent and its children.

\medskip
\noindent\emph{Preparing the sum for telescoping.}
The ordinary averages $f_I$ need not be monotone along nested
intervals. We therefore introduce the maximal average
\[
 \omega_I=\sup_{L\supseteq I}f_L
\]
and use
\[
 f_I^p\leq\omega_I^{p-1}f_I.
\]
As $I$ shrinks, the supremum is taken over a larger collection,
so $\omega_I$ can only increase; it is also bounded by $Mf$ on
$I$. Retaining one ordinary average is important: the quantity
$f_I|I|=\int_If$ is additive over the children and increases when
$I$ is replaced by a containing interval. When $\omega_I$ grows
rapidly from parent to child, geometric summation bounds the
corresponding part of the sparse operator pointwise by $Mf$.
This gives the preliminary removal described above.

For the weight term, we need a form of the tail that can be followed
from one generation to the next. For a dyadic grid $\D$, define
\[
 J_I^{\D}=\sum_{k\geq0}q^k w_{I^{(k)}},
 \qquad q=2^{-(p-1)},
\]
where $I^{(k)}$ is the $k$th dyadic ancestor of $I$.
This discrete analogue of the tail in the $C_p$ condition satisfies
$J^{\D}_{I_\pm}=qJ_I^{\D}+w_{I_\pm}$. This parent--child relation, together
with the monotonicity of $\omega_I$ and the additivity of
$\int_If$, is what makes the telescoping argument possible.

\medskip
\noindent\emph{Faithful intervals and three grids.}
There is a difficulty with this replacement: the $C_p$ condition
controls the $L^r$ average by the full tail, which need not be
bounded by the dyadic tail in a given grid. We call an interval
\emph{faithful} if its full tail is controlled by its dyadic tail.
On faithful intervals, the $C_p$ condition therefore gives the
bound by $J_I^{\D}$ needed for the telescoping argument.

The original sparse family need not consist of faithful intervals,
and we cannot simply discard those that are not faithful.
The three-grid theorem resolves this difficulty. After enlarging
any interval by a fixed factor, we can place it in one of three
fixed grids so that the enlarged interval is faithful and its
dyadic tail also controls the full tail of the original interval.
This construction dominates the original sparse operator by three
sparse operators on faithful intervals, with a fixed loss in the
sparseness constant. It also allows us to transfer terms between
grids later in the proof.

\medskip
\noindent\emph{Telescoping and weight drops.}
On faithful intervals where $J_I^{\D}$ is comparable to $w_I$, the
$C_p$ condition gives the local comparison required in the second
assertion of Lemma~\ref{lem:sparse}, and that assertion bounds
the corresponding sum by $\int_{\R}(Mf)^pw$.
After treating this part and the rapid-growth part separately,
we are left with large dyadic tails and slow growth of $\omega_I$.

For this remaining sum, the telescoping argument transfers the
contributions to parents with small tails but a child with a large
tail, up to terms already controlled by $\int_{\R}(Mf)^pw$.
The parent--child relation for $J_I^{\D}$ forces a fixed relative drop
in the average of $w$ from each such parent to that child.
Thus the transformed sum has additional structure: its intervals
are selected by a drop in the average of the weight.

\medskip
\noindent\emph{Carleson embedding and absorption.}
An entropy argument, based on the convexity of $t\log t$,
controls how often these drops can occur. Where the ratio of
the $L^r$ average of $w$ to its ordinary average is bounded,
it yields a weighted Carleson estimate: the sum of the
$w$-masses of the selected intervals contained in any test
interval is bounded by a constant times its $w$-mass.
The Carleson embedding lemma then controls this part of the
sum by $\int_{\R}(Mf)^pw$.

For intervals where that ratio is large, we use the three-grid
construction again to transfer the remaining terms to enlarged
faithful intervals. The $C_p$ condition supplies a coefficient
proportional to the reciprocal of the ratio, so this coefficient
can be made small. Apart from contributions controlled by the
preceding estimates, the transferred terms have the same form
as the sum we are trying to bound. Summing over the three grids
gives a small multiple of that same total sum, which we move
to the left-hand side. This absorption completes the argument.

\medskip
\noindent\emph{Organization of the paper.}

Section~\ref{sec:prelim} collects the elementary facts used in the proof.
Section~\ref{sec:sparse} reduces the sparse estimate to a sum involving
dyadic tails, and Section~\ref{sec:tree} transfers this sum to the
boundary of the stopping trees. Entropy estimates for these boundary
intervals are proved in Section~\ref{sec:entropy}. Section~\ref{sec:grids}
constructs the assignment to three fixed grids.
Section~\ref{sec:completion} completes the absorption argument and
the proof of Theorem~\ref{thm:main}. Section~\ref{sec:conclusion}
discusses the Fefferman--Stein inequality, the class $SC_p$, and
the extension to higher dimensions and Calder\'on--Zygmund operators.
Appendix~\ref{app:sparse-sum} explains why the full sum in
\eqref{eq:sparse-r} cannot be estimated by the maximal-function integral.

\section{Preliminary facts}\label{sec:prelim}

All intervals are bounded and may be taken half-open. By a
\emph{dyadic grid} we mean a dyadic lattice in the sense of
\cite[Definition~2.1]{LN}. Thus a collection $\D$ of intervals is a dyadic grid
if it satisfies the following three conditions:
\begin{enumerate}
\item both dyadic children of every $I\in\D$ belong to $\D$;
\item any two members of $\D$ have a common dyadic ancestor in $\D$;
\item every compact subset of $\R$ is contained in some $I\in\D$.
\end{enumerate}
In particular, every finite family of grid intervals has a common
ancestor. After a fixed dilation, we take the interval lengths
in the main argument to be $2^k$, $k\in\mathbb Z$.
For $I\in\D$, write $I_-,I_+$ for its children,
$\widehat I$ for its parent, and $I^{(m)}$ for its $m$th ancestor.

We write $w(E)=\int_Ew$ and $w_I=|I|^{-1}w(I)$.
The notation $X\lesssim_{a,b}Y$ means that $X\leq CY$ with a constant
depending only on the displayed parameters. Constants in the main
argument may depend on $p,r,K,\eta$, but never on the minimum scale or on
the cardinality of a sparse family.

It suffices to prove Theorem~\ref{thm:main} for finite sparse families
and nonnegative bounded compactly supported $f$, since
$|A_{\Sfam}f|\leq A_{\Sfam}|f|$. Increasing finite subfamilies and
the functions $\min(|f|,n)\one_{[-n,n]}$ then give the
general result by monotone convergence. We may assume that
\[
 0<\int_{\R}(Mf)^p w<\infty.
\]

Under this assumption, all the tails in the $C_p$ condition are finite.
Indeed, choose an interval $L$ with $f_L>0$. Then
\[
 Mf\geq f_L M\one_L.
\]
On $L$ this is immediate; outside $L$, use the shortest interval
containing $L$ and the point in question. The maximal functions of
the indicators of any two fixed intervals are comparable, with
constants depending on those intervals. It follows that
$\int_{\R}(M\one_I)^pw<\infty$ for every $I$. The $C_p$ condition
therefore also gives $w^r\in L^1_{\mathrm{loc}}$.

We will divide by averages of $w$. To avoid zero denominators, we
reduce to $w>0$ almost everywhere. For $\epsilon>0$, Minkowski's
inequality gives
\[
 \left(\avg{(w+\epsilon)^r}{I}\right)^{1/r}
 \leq \frac K{|I|}\int_{\R}(M\one_I)^p w+\epsilon
 \leq \frac{\max(K,1)}{|I|}
       \int_{\R}(M\one_I)^p(w+\epsilon).
\]
Thus $w+\epsilon$ satisfies the same type of $C_p$ condition, uniformly
in $\epsilon$. We also need the maximal-function integral to remain
finite after replacing $w$ by $w+\epsilon$.
Since $p>1$ and $f$ is bounded and compactly supported, the
unweighted integral $\int_{\R}(Mf)^p$ is finite. Together with
our assumption on the original weighted integral, this gives
\[
 \int_{\R}(Mf)^p(w+\epsilon)
 =\int_{\R}(Mf)^pw+\epsilon\int_{\R}(Mf)^p<\infty.
\]
After proving the estimate for $w+\epsilon$, with a constant
independent of $\epsilon$, we let $\epsilon\downarrow0$.

For every interval $I$, define the maximal average
\[
 \omega_I=\sup_{L\supseteq I}f_L.
\]
The supremum is over all intervals containing $I$. This choice makes
$\omega_I$ independent of any grid, which will allow us to compare
intervals from different grids. Directly from the definition,
\begin{equation}\label{eq:omega-basic}
 f_I\leq\omega_I\leq\inf_{x\in I}Mf(x),
 \qquad
 \omega_{\widehat I}\leq\omega_I\leq2\omega_{\widehat I}.
\end{equation}
For the last inequality, enlarge an interval containing $I$ to include
$\widehat I$; its length increases by at most a factor of two.
More generally, the same argument shows that
\begin{equation}\label{eq:omega-comparable}
 I\subset Q,\quad |Q|\leq a|I|
 \quad\Longrightarrow\quad
 \omega_Q\leq\omega_I\leq a\omega_Q.
\end{equation}
We also have the useful pointwise inequality
\begin{equation}\label{eq:omega-indicator}
 \omega_I M\one_I\leq Mf.
\end{equation}
Indeed, for every $L\supseteq I$,
$Mf\geq f_LM\one_L\geq f_LM\one_I$, and we take the supremum over $L$.
This estimate will control the top intervals of the stopping trees.

Rapid growth of $\omega_I$ provides a separate, elementary source of
summability. We state it for several generations at once because
changing grids can move a jump away from the immediate parent.

\begin{lemma}\label{lem:growth}
Let $m\geq1$ be an integer, and let $a>1$ and $s>0$.
If $\mathscr F\subset\D$ satisfies
$\omega_I\geq a\omega_{I^{(m)}}$ for every $I\in\mathscr F$, then
\[
 \sum_{I\in\mathscr F}\omega_I^s\one_I
 \leq \frac{m}{1-a^{-s}}(Mf)^s
 \qquad\text{almost everywhere}.
\]
\end{lemma}

\begin{proof}
For every $I\in\mathscr F$,
\[
 (1-a^{-s})\omega_I^s
 \leq\omega_I^s-\omega_{I^{(m)}}^s.
\]
All differences on the right are nonnegative by
\eqref{eq:omega-basic}. At a fixed point outside the grid boundaries,
the intervals containing that point form a chain. The sum of these
differences over all scales telescopes in $m$ separate chains of
generations, and is at most $m(Mf)^s$. Finite partial sums give the
same bound, so monotone convergence proves the claim.
\end{proof}

We will also use the Carleson formulation of sparseness. For any
family $\Sfam\subset\D$, finite or infinite, the definition in the
introduction is equivalent to
\begin{equation}\label{eq:sparse-carleson}
 \sum_{\substack{I\in\Sfam\\I\subseteq L}}|I|
 \leq\eta^{-1}|L|
 \qquad(L\in\D).
\end{equation}
See \cite[Lemma~6.3]{LN}. Summing the measures of the disjoint
sets $E_I$ also shows that \eqref{eq:sparse-carleson} holds for
every bounded interval $L$, even when $L\notin\D$.

The weighted version needed below concerns coefficients rather than
lengths. We include the short proof so that all later entropy
estimates can be used through the same elementary principle.

\Needspace{11\baselineskip}
\begin{lemma}\label{lem:embedding}
Suppose $a_I\geq0$, $I\in\D$, and
\[
 \sum_{\substack{I\in\D\\I\subseteq L}}a_I
 \leq\Lambda w(L)
 \qquad(L\in\D).
\]
Then, for every nonnegative measurable function $g$,
\[
 \sum_{I\in\D}a_I\,\essinf_I g
 \leq\Lambda\int_{\R}g\,w.
\]
In particular,
\begin{equation}\label{eq:embedding}
 \sum_{I\in\D}a_I\omega_I^p
 \leq\Lambda\int_{\R}(Mf)^pw.
\end{equation}
\end{lemma}

\begin{proof}
First restrict to finitely many coefficients. For $t>0$, the maximal
intervals with $\essinf_Ig>t$ are disjoint and lie, up to null sets,
in $\{g>t\}$. Applying the hypothesis on these maximal intervals gives
\[
 \sum_{\essinf_Ig>t}a_I\leq\Lambda w(\{g>t\}).
\]
Integration in $t$ proves the first assertion. Increasing finite
subfamilies give the general case, and \eqref{eq:omega-basic} gives
\eqref{eq:embedding}.
\end{proof}

\section{Sparse estimates and dyadic tails}\label{sec:sparse}

The first lemma is the only place where the sparseness of the
original operator is used analytically. It converts the operator
estimate into a sum of local $L^r$ averages of the weight. When these
averages are comparable to $w_I$, the maximal function already
controls the sum.

\begin{lemma}\label{lem:sparse}
Let $\Sfam\subset\D$ be finite and $\eta$-sparse. Then
\begin{equation}\label{eq:sparse-r}
 \int_{\R}(A_{\Sfam}f)^p w
 \lesssim_{p,r,\eta}
 \sum_{I\in\Sfam}f_I^p
       \left(\avg{w^r}{I}\right)^{1/r}|I|.
\end{equation}
If, in addition,
\[
 \left(\avg{w^r}{I}\right)^{1/r}\leq K_0w_I
 \qquad(I\in\Sfam),
\]
then
\[
 \sum_{I\in\Sfam}f_I^p
       \left(\avg{w^r}{I}\right)^{1/r}|I|
 \lesssim_{p,r,\eta,K_0}\int_{\R}(Mf)^p w.
\]
\end{lemma}

\begin{proof}
We use the standard principle that skipping layers improves the
Carleson constant. More precisely, a $\Lambda$-Carleson family can be
partitioned into $m$ subfamilies with Carleson constant at most
$1+(\Lambda-1)/m$; see \cite[Lemma~6.6]{LN}. Taking
$\Lambda=\eta^{-1}$ and $m$ sufficiently large, we may therefore
reduce to a family satisfying
\begin{equation}\label{eq:strong-sparse}
 \sum_{\substack{J\in\Sfam\\J\subsetneq I}}|J|
 \leq\theta|I|
 \qquad(I\in\Sfam),
\end{equation}
where $\theta>0$ is as small as needed. Minkowski's inequality and
summation over the finitely many subfamilies recover the original
family.

Let $\Sfam_k$ be the $k$th generation of this family, starting with its
maximal intervals, and set $U_k=\bigcup_{I\in\Sfam_k}I$.
For $I\in\Sfam_k$, iteration of \eqref{eq:strong-sparse} gives
\[
 |I\cap U_{k+\ell}|\leq\theta^\ell|I|.
\]
Put $E_{I,\ell}=I\cap(U_{k+\ell}\setminus U_{k+\ell+1})$.
For each fixed $\ell$, these sets are pairwise disjoint as $I$
varies, and $I=\bigsqcup_{\ell\geq0}E_{I,\ell}$.
Writing $r'=r/(r-1)$, Minkowski's and H\"older's inequalities give
\begin{align*}
 \|A_{\Sfam}f\|_{L^p(w)}
 &\leq\sum_{\ell\geq0}
       \left(\sum_{I\in\Sfam}f_I^p w(E_{I,\ell})\right)^{1/p}\\
 &\leq\sum_{\ell\geq0}\theta^{\ell/(pr')}
       \left(\sum_{I\in\Sfam}f_I^p
       \left(\avg{w^r}{I}\right)^{1/r}|I|\right)^{1/p}.
\end{align*}
The geometric series proves \eqref{eq:sparse-r}.

For the second assertion, choose the partition so that
$K_0\theta^{1/r'}\leq1/2$. Let $F_I$ be $I$ minus the union of its
immediate $\Sfam$-children. The sets $F_I$ are pairwise disjoint and
\[
 w(I\setminus F_I)
 \leq\theta^{1/r'}\left(\avg{w^r}{I}\right)^{1/r}|I|
 \leq\tfrac12w(I).
\]
Consequently,
\[
 \sum_{I\in\Sfam}f_I^p
       \left(\avg{w^r}{I}\right)^{1/r}|I|
 \leq2K_0\sum_{I\in\Sfam}f_I^p w(F_I)
 \leq2K_0\int_{\R}(Mf)^p w.
\]
Summing over the subfamilies proves the second assertion.
\end{proof}

To use the $C_p$ condition in a dyadic argument, we introduce the
tail obtained by following the ancestors of an interval:
\[
 q=2^{-(p-1)},\qquad
 J_I=J_I^{\D}=\sum_{k\geq0}q^k w_{I^{(k)}}.
\]
The superscript will be omitted when the grid is fixed. These tails
are finite by the preliminary reductions. If $M_{\D}$ denotes the
dyadic maximal operator, then
\begin{equation}\label{eq:J-exact}
 J_I|I|=\frac1{1-2^{-p}}\int_{\R}(M_{\D}\one_I)^p w.
\end{equation}
Indeed, for $x\in I^{(m)}\setminus I^{(m-1)}$ the dyadic maximal
function is $2^{-m}$, and on $I$ it is one. Thus
\[
 \sum_{k\geq0}2^{-kp}\one_{I^{(k)}}
 =\frac{(M_{\D}\one_I)^p}{1-2^{-p}}.
\]
The ancestors of $I$ exhaust $\R$, by our definition of a grid.
Integration proves \eqref{eq:J-exact}.
The definition also gives
\begin{equation}\label{eq:J-recursion}
 J_{I_\pm}=qJ_I+w_{I_\pm},
 \qquad J_I\geq w_I.
\end{equation}

The dyadic tail can miss part of the full tail. We call $I\in\D$
\emph{faithful}, with constant $F$, if
\begin{equation}\label{eq:faithful}
 \frac1{|I|}\int_{\R}(M\one_I)^pw\leq FJ_I^{\D}.
\end{equation}
This is a property of the interval together with its grid.
Section~\ref{sec:grids} will produce faithful intervals with one
constant $F$ depending only on $p$.

Fix a constant $c>1$, depending only on $p$, so large that
\begin{equation}\label{eq:choice-c}
 q+\frac2c<(2/3)^{p-1}.
\end{equation}
Such a choice is possible because $q=(1/2)^{p-1}<(2/3)^{p-1}$.
The significance of this inequality will become clear in the
stopping-tree calculation.

There are two easy parts of the sparse estimate. On faithful intervals
with $J_I\leq cw_I$, the $C_p$ condition gives the local comparison
required in Lemma~\ref{lem:sparse}. Intervals with
$\omega_I\geq\frac32\omega_{\widehat I}$ are controlled by
Lemma~\ref{lem:growth}. We are therefore left with large tails and
slow growth of $\omega_I$.

On the remaining faithful intervals, Lemma~\ref{lem:sparse} and the
$C_p$ condition lead to terms of the form $f_I^pJ_I|I|$.
Before defining the sum that we will estimate, we deliberately
enlarge these terms using
\[
 f_I^p\leq\omega_I^{p-1}f_I.
\]
The maximal average $\omega_I$ has the useful monotonicity and
growth properties established in Section~\ref{sec:prelim}.
Retaining one ordinary average $f_I$, however, preserves the
quantity $f_I|I|=\int_If$. Its additivity over the two children
will be used together with the recursion for $J_I$ in the
stopping-tree calculation. It also satisfies
\[
 f_I|I|=\int_If\leq\int_Qf=f_Q|Q|
 \qquad(I\subset Q),
\]
which will allow us to transfer a contribution to a larger
interval in another grid. Replacing $f_I^p$ by $\omega_I^p$
would lose this direct comparison of the integrals of $f$.

With this choice in mind, fix a minimum interval length $h=2^k$
and define
\begin{equation}\label{eq:A}
 \mathcal A_h^{\D}
 =
 \sum_{\substack{I\in\D,\ |I|\geq h\\
                  J_I>cw_I,\ \omega_I<\frac32\omega_{\widehat I}}}
       \omega_I^{p-1}f_IJ_I|I|.
\end{equation}
Only the summation is truncated: each $J_I$ still includes all
ancestors.

This sum is finite. To check the upper scales, let $I_*$ contain
the support of $f$. If $f_I>0$, the convex hull of $I$ and $I_*$
has length at most $|I|+|I_*|$, so
\[
 \omega_I\geq\frac{\|f\|_1}{|I|+|I_*|},
 \qquad
 \omega_{\widehat I}\leq\frac{\|f\|_1}{2|I|}.
\]
Consequently, $\omega_I<\frac32\omega_{\widehat I}$ implies
$|I|<3|I_*|$. There are only finitely many scales between $h$ and
this upper bound, and only finitely many intervals with $f_I>0$
at each such scale. This observation will justify the absorption
at the end of the proof.

\begin{proposition}\label{prop:sparse-reduction}
Let $\Sfam\subset\D$ be a finite $\eta$-sparse family of faithful
intervals with constant $F$, and suppose $|I|\geq h$ for
$I\in\Sfam$. Then
\[
 \int_{\R}(A_{\Sfam}f)^p w
 \lesssim_{p,r,K,\eta,F}
 \int_{\R}(Mf)^p w+\mathcal A_h^{\D}.
\]
\end{proposition}

\begin{proof}
For the subfamily on which $J_I\leq cw_I$, the $C_p$ condition
and faithfulness \eqref{eq:faithful} give
\[
 \left(\avg{w^r}{I}\right)^{1/r}\leq KFJ_I\leq KFcw_I.
\]
Both assertions of Lemma~\ref{lem:sparse} bound this part of the
operator by the required maximal-function integral.

For the subfamily $\Sfam'$ with $J_I>cw_I$ and
$\omega_I<\frac32\omega_{\widehat I}$, the same inequalities and
\eqref{eq:sparse-r} give
\[
 \int_{\R}(A_{\Sfam'}f)^p w
 \lesssim\sum_{I\in\Sfam'}f_I^pJ_I|I|
 \leq\mathcal A_h^{\D}.
\]
Here we used $f_I\leq\omega_I$.

On the remaining subfamily $\Sfam''$,
$\omega_I\geq\frac32\omega_{\widehat I}$.
Lemma~\ref{lem:growth}, with $m=s=1$ and $a=3/2$, gives
\[
 A_{\Sfam''}f
 \leq\sum_{I\in\Sfam''}\omega_I\one_I
 \leq3Mf.
\]
Combining the three parts proves the proposition.
\end{proof}

The sum $\mathcal A_h^{\D}$ runs over all intervals satisfying its
two conditions, not just a sparse family. From this point until the
final transfer of the original operator, we work with the full
dyadic trees.

\section{Stopping at small dyadic tails}\label{sec:tree}

We now estimate $\mathcal A_h^{\D}$ by summing parent--child
differences over the intervals with large tails, $J_I>cw_I$.
When both an interval and its parent have large tails, the
corresponding terms cancel. Apart from a possible top-interval
contribution, the positive terms left after cancellation come
from intervals $I$ with large tails whose parents $P$ have small
tails, $J_P\leq cw_P$. We estimate these terms by a sum over
their parents $P$; the other contributions are controlled by
$\int_{\R}(Mf)^pw$.

For a fixed grid and minimum interval length $h$, let
\begin{equation}\label{eq:boundary-family}
 \Pfamily_h^{\D}
 =
 \left\{P\in\D:\ |P|\geq2h,\ J_P\leq cw_P,\
       J_I>cw_I\ \text{for some child }I\text{ of }P\right\}.
\end{equation}
These are the parents with a small tail and at least one child
with a large tail. We call them \emph{boundary parents}: each
marks a transition from a small tail at $P$ to a large tail at
one of its children as we descend through the dyadic grid.

The estimate below retains the product $\omega_P^{p-1}f_P$
introduced in Section~\ref{sec:sparse}.

\begin{proposition}\label{prop:tree}
For every dyadic grid $\D$ and every minimum length $h=2^k$,
\begin{equation}\label{eq:tree}
 \mathcal A_h^{\D}
 \leq C_p\int_{\R}(Mf)^p w
     +C_p\sum_{P\in\Pfamily_h^{\D}}
               \omega_P^{p-1}f_Pw(P).
\end{equation}
Moreover, every $P\in\Pfamily_h^{\D}$ satisfies
\begin{equation}\label{eq:weight-drop}
 \min(w_{P_-},w_{P_+})\leq\tau w_P,
 \qquad
 \tau=\frac{qc}{c-1}<1.
\end{equation}
\end{proposition}

\begin{proof}
Choose a single interval $L\in\D$, with $|L|\geq h$, containing
every interval that contributes a nonzero summand to
$\mathcal A_h^{\D}$. This is possible because there are only
finitely many such intervals, as shown after \eqref{eq:A}, and
a finite family in $\D$ has a common ancestor.

Inside $L$, consider the finite family
\[
 \mathscr B
 =\{I\in\D:I\subseteq L,\ |I|\geq h,\ J_I>cw_I\}.
\]
The interval $L$ only specifies the region in which we work;
it is not automatically included in $\mathscr B$. In particular,
if $J_L\leq cw_L$, then $L\notin\mathscr B$, even if
$\omega_L<\frac32\omega_{\widehat L}$.
Thus $\mathcal A_h^{\D}$ is the sum over those $I\in\mathscr B$
for which $\omega_I<\frac32\omega_{\widehat I}$.
The family $\mathscr B$ also includes intervals where this inequality
fails. We will sum parent--child differences over the whole family
$\mathscr B$ to obtain cancellation, so both types of intervals
must be estimated.
For every grid interval $I$, put
\[
 m_I=\omega_{\widehat I}^{p-1}f_IJ_I|I|.
\]
We will compare $m_I$ with the sum of its two child terms.
The use of $\omega_{\widehat I}$ makes their common factor equal
to $\omega_I^{p-1}$.

Separate the following three contributions:
\[
 \begin{aligned}
 A&=m_L,\\
 B&=\sum_{\substack{I\in\mathscr B,\ I\ne L\\
                    \widehat I\notin\mathscr B}}m_I,\\
 C&=2\sum_{\substack{I\in\mathscr B\\
                     \omega_I\geq\frac32\omega_{\widehat I}}}
           \omega_I^{p-1}f_Iw(I).
 \end{aligned}
\]
Here $A$ bounds the possible contribution of the top interval.
If $L\notin\mathscr B$, this contribution is absent, and the
nonnegative term $A=m_L$ merely enlarges the upper bound below.
The sum $B$ collects large-tail intervals whose parents have
small tails, and $C$ accounts for intervals where $\omega$
grows rapidly. We first show that
\begin{equation}\label{eq:three-contributions}
 \kappa\mathcal A_h^{\D}\leq A+B+C,
 \qquad
 \kappa=(2/3)^{p-1}-q-\frac2c>0,
\end{equation}
and then estimate $A$, $B$, and $C$ separately. Positivity of
$\kappa$ follows from \eqref{eq:choice-c}.

\medskip
\noindent\emph{Reduction to the three contributions.}
Using \eqref{eq:J-recursion} and the identities for child
averages, we obtain
\begin{align*}
 m_{I_-}+m_{I_+}
 &=\omega_I^{p-1}|I|
   \left(qJ_If_I+
      \frac{f_{I_-}w_{I_-}+f_{I_+}w_{I_+}}2\right)\\
 &\leq\omega_I^{p-1}f_I|I|(qJ_I+2w_I).
\end{align*}
The last step uses nonnegativity together with
$f_{I_-}+f_{I_+}=2f_I$ and $w_{I_-}+w_{I_+}=2w_I$.
For $I\in\mathscr B$ with
$\omega_I<\frac32\omega_{\widehat I}$, namely the intervals
retained in $\mathcal A_h^{\D}$, this gives
\begin{equation}\label{eq:tree-local}
 m_I-m_{I_-}-m_{I_+}
 \geq\kappa\omega_I^{p-1}f_IJ_I|I|.
\end{equation}
We next consider the intervals $I\in\mathscr B$ with
$\omega_I\geq\frac32\omega_{\widehat I}$. These do not enter
$\mathcal A_h^{\D}$, but they occur in the sum over
$\mathscr B$ used for cancellation. For these intervals,
\eqref{eq:omega-basic} gives
$\omega_{\widehat I}^{p-1}\geq q\omega_I^{p-1}$, and hence
\begin{equation}\label{eq:tree-loss}
 m_I-m_{I_-}-m_{I_+}
 \geq-2\omega_I^{p-1}f_Iw(I).
\end{equation}
Summing \eqref{eq:tree-local} over the slow-growth intervals in
$\mathscr B$ and \eqref{eq:tree-loss} over the rapid-growth
intervals in $\mathscr B$ yields
\[
 \kappa\mathcal A_h^{\D}
 \leq\sum_{I\in\mathscr B}(m_I-m_{I_-}-m_{I_+})+C.
\]

To see the cancellation, suppose that both $I$ and $\widehat I$
belong to $\mathscr B$. The summand indexed by $I$ contains
$+m_I$, while the summand indexed by $\widehat I$ contains
$-m_I$, because $I$ is a child of $\widehat I$. These two
occurrences of the nonnegative number $m_I$ cancel.
The only possible uncancelled terms with a plus sign are
$m_L$, present only if $L\in\mathscr B$, and the terms defining
$B$. Every other uncancelled term
has the form $-m_J\leq0$ and can be discarded in an upper bound.
Therefore
\[
 \sum_{I\in\mathscr B}(m_I-m_{I_-}-m_{I_+})
 \leq m_L+
       \sum_{\substack{I\in\mathscr B,\ I\ne L\\
                        \widehat I\notin\mathscr B}}m_I
 =A+B.
\]
This proves \eqref{eq:three-contributions}. The argument includes
the smallest retained intervals: their children have length
$h/2$ and occur only among the discarded negative terms.

\medskip
\noindent\emph{Estimate of $A$.}
By \eqref{eq:omega-basic}, \eqref{eq:J-exact}, and
\eqref{eq:omega-indicator},
\[
 A=\omega_{\widehat L}^{p-1}f_LJ_L|L|
 \leq\omega_L^pJ_L|L|
 \leq\frac1{1-2^{-p}}\int_{\R}(Mf)^pw.
\]

\medskip
\noindent\emph{Estimate of $B$.}
If $I$ occurs in $B$, its parent $P=\widehat I$ lies inside $L$
and has length at least $2h$. Since $P\notin\mathscr B$, we have
$J_P\leq cw_P$. Thus $P\in\Pfamily_h^{\D}$.
Using $w_I\leq2w_P$ and the tail recursion, we obtain
\[
 J_I=qJ_P+w_I\leq(qc+2)w_P.
\]
Since $f_I|I|\leq f_P|P|$, it follows that
\begin{equation}\label{eq:entrance-cost}
 m_I\leq(qc+2)\omega_P^{p-1}f_Pw(P).
\end{equation}
Each parent has at most two children. Summing
\eqref{eq:entrance-cost} therefore gives
\[
 B\leq2(qc+2)\sum_{P\in\Pfamily_h^{\D}}
                     \omega_P^{p-1}f_Pw(P).
\]
This is the boundary sum in the statement, with the ordinary
average $f_P$ still present.

\medskip
\noindent\emph{Estimate of $C$.}
Use $f_I\leq\omega_I$ and apply Lemma~\ref{lem:growth} with
$m=1$, $a=3/2$, and $s=p$. We get
\[
 C\leq2\sum_{\substack{I\in\D\\
              \omega_I\geq\frac32\omega_{\widehat I}}}
                 \omega_I^pw(I)
 \lesssim_p\int_{\R}(Mf)^pw.
\]
Substituting the three estimates into
\eqref{eq:three-contributions} proves \eqref{eq:tree}.

\medskip
\noindent\emph{The weight drop at a boundary parent.}
It remains to prove \eqref{eq:weight-drop}. If $I$ is a child
of $P\in\Pfamily_h^{\D}$ with $J_I>cw_I$, then
\[
 cw_I<J_I=qJ_P+w_I\leq qc\,w_P+w_I.
\]
Rearranging gives $w_I<\frac{qc}{c-1}w_P=\tau w_P$.
Finally, \eqref{eq:choice-c} implies $q+2/c<1$, so
$c(1-q)>1$ and hence $\tau<1$.
\end{proof}

The boundary intervals in \eqref{eq:tree} need not be sparse with
respect to Lebesgue measure. Their useful structure is instead
\eqref{eq:weight-drop}: each one witnesses a fixed relative change
in the average of the weight. The next section turns that change
into a weighted counting estimate.

\section{Entropy and weight drops}\label{sec:entropy}

For an interval $I$, let
\[
 R_I=\frac{\left(\avg{w^r}{I}\right)^{1/r}}{w_I}.
\]
By H\"older's inequality, $R_I\geq1$. This ratio measures how much
larger the local $L^r$ average is than the ordinary average. The
entropy estimate below is useful when this ratio is bounded.
When it is large, the $C_p$ condition will instead force a large
tail, and we will use that in the absorption argument.

For $0<\tau<1$ and a dyadic grid $\D$, define
\[
 \Ofamily_\tau^{\D}
 =\{I\in\D:\min(w_{I_-},w_{I_+})\leq\tau w_I\}.
\]
Thus $\Ofamily_\tau^{\D}$ consists of intervals with a definite drop
in the average of $w$ on one child. By \eqref{eq:weight-drop},
the boundary family $\Pfamily_h^{\D}$ is contained in this family
for the fixed value of $\tau$ chosen in Section~\ref{sec:tree}.

The containing interval in the next lemma is deliberately allowed
to lie outside $\D$. This will be essential when an interval is
moved from one grid to another.

\begin{lemma}\label{lem:entropy}
Let $0<\tau<1$, let $\D$ be a dyadic grid, and let $L$ be any
bounded interval. Then
\begin{equation}\label{eq:entropy}
 \sum_{\substack{I\in\Ofamily_\tau^{\D}\\I\subseteq L}}w(I)
 \leq C_\tau\int_Lw\log\frac{w}{w_L}
 \leq C_\tau r'(\log R_L)\,w(L),
 \qquad r'=\frac r{r-1}.
\end{equation}
\end{lemma}

\begin{proof}
Put $\phi(t)=t\log t$ for $t>0$, with $\phi(0)=0$.
For $I\in\D$, the identity
$w_I=(w_{I_-}+w_{I_+})/2$ and convexity give
\[
 \phi(w_I)\leq\frac{\phi(w_{I_-})+\phi(w_{I_+})}{2}.
\]
Set
\[
 d_I=\frac{|I|}{2}
       \bigl(\phi(w_{I_-})+\phi(w_{I_+})-2\phi(w_I)\bigr)\geq0.
\]
Strict convexity shows that $d_I=0$ exactly when the two child
averages are equal. A fixed relative drop on one child gives the
following quantitative lower bound.
If $I\in\Ofamily_\tau^{\D}$, write the two child averages as
$tw_I$ and $(2-t)w_I$, with $0\leq t\leq\tau$. Then
\[
 d_I=\frac{w(I)}2
       \bigl(t\log t+(2-t)\log(2-t)\bigr)
 \geq\kappa_\tau w(I),
\]
where
\[
 \kappa_\tau
 =\tfrac12\bigl(\tau\log\tau+(2-\tau)\log(2-\tau)\bigr)>0.
\]
The lower bound follows because the expression in parentheses
decreases on $[0,1]$.

Suppose first that $L\in\D$, and fix $n\geq1$. Sum $d_I$ over
the dyadic subintervals of $L$ in generations $0,\ldots,n-1$.
The child terms from one generation cancel the corresponding
parent terms from the next, leaving
\[
 \sum_{\substack{I\in\D,\ I\subseteq L\\|I|>2^{-n}|L|}}d_I
 =\sum_{\substack{J\in\D,\ J\subseteq L\\|J|=2^{-n}|L|}}
       |J|\phi(w_J)-|L|\phi(w_L).
\]
The intervals $J$ on the right form the bottom generation and
partition $L$. On each such interval, Jensen's inequality for
the convex function $\phi$ gives
\[
 \phi(w_J)
 =\phi\left(\frac1{|J|}\int_Jw\right)
 \leq\frac1{|J|}\int_J\phi(w).
\]
Multiplying by $|J|$ and summing over these disjoint intervals
bounds the sum over $J$ by $\int_L\phi(w)$. Therefore
\[
 \sum_{\substack{I\in\D,\ I\subseteq L\\|I|>2^{-n}|L|}}d_I
 \leq\int_L\phi(w)-|L|\phi(w_L)
 =\int_Lw\log\frac{w}{w_L}.
\]
Since $d_I\geq0$ for every $I$ and
$d_I\geq\kappa_\tau w(I)$ for $I\in\Ofamily_\tau^{\D}$,
we may discard the other intervals and let $n\to\infty$.
This proves the first inequality of \eqref{eq:entropy} for
dyadic $L$, with $C_\tau=\kappa_\tau^{-1}$.

For an arbitrary $L$, first take a finite subfamily of the sum
on the left and let $L_1,\ldots,L_N$ be its maximal intervals.
They are disjoint dyadic intervals contained in $L$. The preceding
argument applies on each $L_a$. Moreover,
\begin{equation}\label{eq:entropy-partition}
 \sum_{a=1}^N\int_{L_a}w\log\frac{w}{w_{L_a}}
 \leq\int_Lw\log\frac{w}{w_L}.
\end{equation}
To verify this, let $E=L\setminus\bigcup_aL_a$. The difference
between the right and left sides is
\[
 \sum_{a=1}^N|L_a|\phi(w_{L_a})
       +\int_E\phi(w)-|L|\phi(w_L).
\]
Suppose first that $|E|>0$, and put
$w_E=|E|^{-1}\int_Ew$. Jensen's inequality on $E$ gives
\[
 \int_E\phi(w)\geq |E|\phi(w_E).
\]
Next, since $L_1,\ldots,L_N,E$ form a partition of $L$, we have
\[
 w_L=\sum_{a=1}^N\frac{|L_a|}{|L|}w_{L_a}
              +\frac{|E|}{|L|}w_E.
\]
The coefficients in this expression are nonnegative and sum to
one. Applying convexity of $\phi$ to this weighted average and
multiplying by $|L|$, we obtain
\[
 |L|\phi(w_L)
 \leq\sum_{a=1}^N|L_a|\phi(w_{L_a})+|E|\phi(w_E)
 \leq\sum_{a=1}^N|L_a|\phi(w_{L_a})+\int_E\phi(w).
\]
Thus the displayed difference is nonnegative, as required.
If $|E|=0$, the coefficients $|L_a|/|L|$ already sum to one,
and the same convexity argument applies with the $E$ term omitted.
This proves \eqref{eq:entropy-partition}. Increasing finite
subfamilies proves the first inequality for every $L$.

For the second inequality, use the normalized measure
\[
 d\mu(x)=\frac{w(x)}{w(L)}\,dx\quad\text{on }L,
 \qquad \mu(L)=1.
\]
Jensen's inequality for the concave function $\log$ says that
$\int_L\log g\,d\mu\leq\log\int_L g\,d\mu$.
Apply it to $g=(w/w_L)^{r-1}$. This choice of power makes
$gw=w^r/w_L^{r-1}$, so that the resulting integral is controlled
by the $L^r$ average of $w$. Using
$\log(w/w_L)=(r-1)^{-1}\log g$, we obtain
\begin{align*}
 \frac1{w(L)}\int_Lw\log\frac{w}{w_L}
 &=\frac1{r-1}\int_L\log g\,d\mu\\
 &\leq\frac1{r-1}
       \log\left(\frac1{w(L)}
                  \int_Lw\left(\frac{w}{w_L}\right)^{r-1}\right)\\
 &=\frac1{r-1}
       \log\left(\frac{|L|^{-1}\int_Lw^r}{w_L^r}\right)\\
 &=r'\log R_L.
\end{align*}
This completes the proof.
\end{proof}

The first consequence controls weight drops wherever the
concentration ratio is bounded. The intervals involved need not
form a sparse family.

\begin{corollary}\label{cor:bounded-R}
For every $T>1$,
\[
 \sum_{\substack{I\in\Ofamily_\tau^{\D}\\R_I\leq T}}
       \omega_I^pw(I)
 \leq C_{\tau,r}\log T\int_{\R}(Mf)^pw.
\]
\end{corollary}

\begin{proof}
Consider a finite subfamily and a dyadic test interval $L$. Choose
the maximal intervals of the subfamily contained in $L$. On each
such interval $Q$, \eqref{eq:entropy} bounds the sum of the weights
of all its descendants in the subfamily by
$C_\tau r'\log T\,w(Q)$. The maximal intervals are disjoint, so
the coefficients $w(I)$ satisfy the hypothesis of
Lemma~\ref{lem:embedding} with
$\Lambda=C_\tau r'\log T$. Applying \eqref{eq:embedding} and
increasing the finite subfamily proves the claim.
\end{proof}

Corollary~\ref{cor:bounded-R} applies when $R_P$ is bounded on
the intervals $P$ where the weight drops. Later we will encounter
intervals for which $R_P$ may be large, but for which we can
choose a containing interval $Q(P)$ of comparable length with
$R_{Q(P)}\leq T$. The next lemma shows that this bound on the
larger interval still suffices to control the sum over the
original intervals $P$.

The intervals $P$ may belong to several grids
$\D_1,\ldots,\D_N$. We require all the assigned intervals $Q(P)$
to belong to one grid $\D_*$. In the proof, we sum the $w$-masses
$w(P)$ of all intervals assigned to the same $Q$ and apply the
Carleson embedding lemma in $\D_*$. The required Carleson estimate
comes from Lemma~\ref{lem:entropy}: it can be applied with an
assigned interval as the containing interval, even when that
interval does not belong to the grid of the original intervals $P$.

\Needspace{14\baselineskip}
\begin{lemma}\label{lem:entropy-transfer}
Let $\D_1,\ldots,\D_N,\D_*$ be dyadic grids. For each $\nu$, let
$\mathscr F_\nu\subset\Ofamily_\tau^{\D_\nu}$
and assign to every $P\in\mathscr F_\nu$ an interval $Q(P)\in\D_*$
such that
\[
 P\subset Q(P),\qquad |Q(P)|\leq a|P|,
 \qquad R_{Q(P)}\leq T,
\]
where $a\geq1$ and $T>1$. Then
\[
 \sum_{\nu=1}^N\sum_{P\in\mathscr F_\nu}\omega_P^pw(P)
 \leq C_{\tau,r}Na^p\log T\int_{\R}(Mf)^pw.
\]
\end{lemma}

\begin{proof}
First assume that each $\mathscr F_\nu$ is finite. For every
$Q\in\D_*$, sum the $w$-masses of the original intervals
assigned to $Q$ by setting
\[
 a_Q=\sum_{\nu=1}^N
           \sum_{\substack{P\in\mathscr F_\nu\\Q(P)=Q}}w(P).
\]
An empty sum is understood to be zero. We claim that these
coefficients are $w$-Carleson in $\D_*$.
Fix a test interval $L\in\D_*$. Among the assigned intervals
$Q(P)\subseteq L$, choose those that are maximal with respect
to inclusion, and denote them by $Q_0$.
If $Q(P)\subseteq Q_0$, then $P\subset Q(P)\subseteq Q_0$.
We may therefore apply Lemma~\ref{lem:entropy} in each grid
$\D_\nu$, using $Q_0$ as the containing interval. This gives
\[
 \sum_{\substack{Q\in\D_*\\Q\subseteq Q_0}}a_Q
 \leq C_\tau N\int_{Q_0}w\log\frac{w}{w_{Q_0}}
 \leq C_\tau Nr'\log T\,w(Q_0).
\]
Here $R_{Q_0}\leq T$ because $Q_0$ is one of the assigned
intervals $Q(P)$. The maximal intervals $Q_0$ are disjoint,
and every $Q\subseteq L$ with $a_Q\ne0$ is contained in one of
them. Hence
\[
 \sum_{Q\subseteq L}a_Q
 \leq C_\tau Nr'\log T\,w(L).
\]
Using \eqref{eq:omega-comparable} and
Lemma~\ref{lem:embedding}, we obtain
\[
 \sum_{\nu=1}^N\sum_{P\in\mathscr F_\nu}\omega_P^pw(P)
 \leq a^p\sum_Q\omega_Q^pa_Q
 \leq C_\tau Nr'a^p\log T\int_{\R}(Mf)^pw.
\]
Increasing finite subfamilies of each $\mathscr F_\nu$ completes
the proof.
\end{proof}

\section{Three fixed dyadic grids}\label{sec:grids}

We now address the difference between a dyadic tail and the full
tail in the $C_p$ condition. A single fixed grid cannot always
capture that tail. Three adjacent grids suffice after enlarging
the interval by a fixed factor. The point is to use the adjacency
at two scales: the scale of the enlarged interval and a larger
scale reaching a given point of the line.

Apply the Three Lattice Theorem \cite[Theorem~3.1]{LN} to an
auxiliary grid $\D_*$ with interval lengths $2^k/3$.
It gives three grids $\D^0,\D^1,\D^2$ such that
$$
 \D^0\cup\D^1\cup\D^2
 =\{[a,a+3\ell):[a,a+\ell)\in\D_*\}.
$$
Each of these grids satisfies the definition in
Section~\ref{sec:prelim}, and its interval lengths are $2^k$.
At a fixed scale $2^k$, the left endpoints of the intervals
from all three grids, taken together, form an arithmetic
progression with spacing $2^k/3$. Within each individual grid,
successive left endpoints are $2^k$ apart, so each grid
contributes every third point of this progression. Thus, as we
move along the line, the endpoints from the three grids occur
in a repeating order, with spacing $2^k/3$ between consecutive
points in the combined set.
In particular, an interval of length less than $2^k/3$ is
contained in a length-$2^k$ interval of at least two grids.

\begin{lemma}\label{lem:selection}
There are constants $C_p,F_p$ such that every interval $P$ of
length $2^\ell$ can be assigned an index $j(P)\in\{0,1,2\}$
and an interval $Q(P)\in\D^{j(P)}$ satisfying
\begin{equation}\label{eq:selection-size}
 P\subset Q(P),\qquad |Q(P)|=8|P|,
\end{equation}
\begin{equation}\label{eq:selection-tail}
 \frac1{|P|}\int_{\R}(M\one_P)^pw
 \leq C_pJ_{Q(P)}^{\D^{j(P)}},
\end{equation}
and
\begin{equation}\label{eq:selection-faithful}
 \frac1{|Q(P)|}\int_{\R}(M\one_{Q(P)})^pw
 \leq F_pJ_{Q(P)}^{\D^{j(P)}}.
\end{equation}
Thus every assigned interval $Q(P)$ is faithful with the same constant
$F_p$. The assignment depends only on $P$, $w$, and the three
fixed grids.
\end{lemma}

\begin{proof}
Set $H=8|P|$. At this scale, each of the three grids partitions
$\R$ into intervals of length $H$. Recall that the endpoints
of these three partitions, taken together, are spaced $H/3$
apart. Since $|P|=H/8<H/3$, the interval $P$ can cross an endpoint
of at most one of the three partitions. Thus, in at least two
grids, $P$ lies entirely inside one interval of length $H$.
For every index $j$ for which such an interval exists, denote
it by $Q_j$. These are the candidates from which we will choose
$Q(P)$.

Fix $x\in\R$ and put
$d=|P|+\operatorname{dist}(x,P)$. The convex hull of $P$ and $x$
has length $d$. To estimate a dyadic maximal function at $x$,
we look for an interval containing both $x$ and one of the
candidates $Q_j$, with length comparable to $d$.
Choose a dyadic length $L$ such that
\[
 \max(H,4d)\leq L<2\max(H,4d)\leq16d.
\]
The convex hull has length $d\leq L/4<L/3$. The same argument
using the spacing of endpoints at scale $L$ shows that at least two grids
have a length-$L$ interval containing both $P$ and $x$.
There are only three grids, so at least one of these also has
a candidate $Q_j$ at scale $H$.

Fix such an index $j$. The length-$L$ interval and $Q_j$ belong
to the same grid and both contain $P$. Since $L\geq H$, the
nesting property implies that the length-$L$ interval contains
$Q_j$. It also contains $x$, so averaging $\one_{Q_j}$ over
this interval gives
\[
 M_{\D^j}\one_{Q_j}(x)
 \geq\frac{H}{L}
 \geq\frac12\,\frac{|P|}{|P|+\operatorname{dist}(x,P)}
 =\frac12M\one_P(x).
\]
The index $j$ may depend on $x$, but the candidates $Q_j$ were
fixed before choosing $x$. Raising to the $p$th power and summing
over all candidates therefore gives a pointwise bound by a fixed
sum of at most three dyadic maximal functions to the $p$th power.
Integrating and using \eqref{eq:J-exact}, we obtain
\[
 \frac1{|P|}\int_{\R}(M\one_P)^pw
 \leq8(2^p-1)\sum_jJ_{Q_j}^{\D^j}.
\]
Choose a candidate maximizing $J_{Q_j}^{\D^j}$, resolving ties by
the smallest index, and call it $Q(P)$. This proves
\eqref{eq:selection-size} and \eqref{eq:selection-tail}.

For nested intervals with $|Q(P)|=8|P|$, the explicit formula
for the maximal function of an indicator gives
\[
 M\one_{Q(P)}\leq8M\one_P.
\]
Consequently,
\[
 \frac1{|Q(P)|}\int_{\R}(M\one_{Q(P)})^pw
 \leq\frac{8^{p-1}}{|P|}\int_{\R}(M\one_P)^pw.
\]
Together with \eqref{eq:selection-tail}, this proves
\eqref{eq:selection-faithful}.
\end{proof}

We henceforth fix this assignment. Two of its consequences explain
how it will be used. The $C_p$ condition and
\eqref{eq:selection-tail} imply
\begin{equation}\label{eq:selection-R}
 R_Pw_P
 =\left(\avg{w^r}{P}\right)^{1/r}
 \leq C_pKJ_{Q(P)}^{\D^{j(P)}}.
\end{equation}
Thus a large concentration ratio $R_P$ forces the selected
dyadic tail to be large compared with $w_P$. Also,
\eqref{eq:selection-size} and \eqref{eq:omega-comparable} give
\begin{equation}\label{eq:selection-f}
 \omega_P\leq8\omega_{Q(P)}.
\end{equation}

The assignment can be applied to the original sparse family as
well as to the boundary intervals. The following proposition
ensures that the resulting operators remain sparse.

\begin{proposition}\label{prop:sparse-transfer}
Let $\Sfam$ be a finite $\eta$-sparse family in any dyadic grid.
For $j=0,1,2$, let
\[
 \Sfam_j=\{Q(I):I\in\Sfam,\ j(I)=j\},
\]
with repeated intervals $Q(I)$ counted only once. Then $\Sfam_j\subset\D^j$
is $\eta/8$-sparse, every interval in $\Sfam_j$ is faithful with
constant $F_p$, and
\begin{equation}\label{eq:sparse-transfer}
 A_{\Sfam}f\leq64\sum_{j=0}^2A_{\Sfam_j}f.
\end{equation}
\end{proposition}

\begin{proof}
Faithfulness follows from \eqref{eq:selection-faithful}.
For each $I\in\Sfam$, the size and containment relations give
\[
 f_I\one_I\leq8f_{Q(I)}\one_{Q(I)}.
\]
For a fixed $Q$, at most eight intervals $I$ in the original
grid satisfy $Q(I)=Q$: they all lie in $Q$, have length $|Q|/8$,
and are disjoint.
This proves \eqref{eq:sparse-transfer}.

For $L\in\D^j$ and each $Q\in\Sfam_j$ contained in $L$, choose
one original interval $I(Q)\in\Sfam$ assigned to $Q$.
These intervals $I(Q)$ are distinct and contained in $L$.
If $E_I$ are the disjoint sets witnessing the
original sparseness, then
\[
 \sum_{\substack{Q\in\Sfam_j\\Q\subseteq L}}|Q|
 =8\sum_{Q\subseteq L}|I(Q)|
 \leq8\eta^{-1}\sum_{Q\subseteq L}|E_{I(Q)}|
 \leq8\eta^{-1}|L|.
\]
The equivalence \eqref{eq:sparse-carleson} proves the asserted
sparseness.
\end{proof}

One more geometric observation is needed. Rapid growth of $\omega$
from the parent of $Q(P)$ to $Q(P)$ need not imply rapid growth
from the parent of $P$ to $P$. It does, however, imply a definite
increase from the sixth ancestor of $P$ to $P$ in the original
grid, which is enough for Lemma~\ref{lem:growth}.

\begin{lemma}\label{lem:transferred-growth}
Let $\D$ be a dyadic grid. Suppose $\mathscr F\subset\D$ consists
of intervals $P$ whose assigned intervals $Q(P)$ satisfy
\[
 \omega_{Q(P)}\geq\tfrac32\omega_{\widehat{Q(P)}},
\]
where the parent of $Q(P)$ is taken in $\D^{j(P)}$. Then
\[
 \omega_{P^{(6)}}\leq\tfrac56\omega_P
 \qquad(P\in\mathscr F),
\]
and
\[
 \sum_{P\in\mathscr F}\omega_P^pw(P)
 \lesssim_p\int_{\R}(Mf)^pw.
\]
\end{lemma}

\begin{proof}
The intervals $P^{(6)}$ and $\widehat{Q(P)}$ both contain $P$,
and have lengths $64|P|$ and $16|P|$, respectively.
Fix any interval $L\supseteq P^{(6)}$, and let $\widetilde L$
be the smallest interval containing both $L$ and
$\widehat{Q(P)}$. These two intervals intersect, since both
contain $P$. Hence there is no gap between them, and
\[
 |\widetilde L|\leq |L|+|\widehat{Q(P)}|
 =|L|+16|P|\leq\tfrac54|L|,
\]
where the last inequality uses $|L|\geq64|P|$.

Since $f\geq0$ and $L\subseteq\widetilde L$, comparison of the
integrals gives
\[
 f_L\leq\frac{|\widetilde L|}{|L|}f_{\widetilde L}
      \leq\tfrac54\omega_{\widehat{Q(P)}}.
\]
Here $f_{\widetilde L}\leq\omega_{\widehat{Q(P)}}$ because
$\widetilde L$ contains $\widehat{Q(P)}$, so its average is
included in the supremum defining this maximal average.
Taking the supremum over all $L\supseteq P^{(6)}$ gives the
first inequality below. The next two follow from the hypothesis
of the lemma and the inclusion $P\subset Q(P)$:
\[
 \omega_{P^{(6)}}
 \leq\tfrac54\omega_{\widehat{Q(P)}}
 \leq\tfrac56\omega_{Q(P)}
 \leq\tfrac56\omega_P,
\]
which proves the first assertion. Apply
Lemma~\ref{lem:growth} with $m=6$, $a=6/5$, and $s=p$, and
integrate against $w$ to obtain the second assertion.
\end{proof}

\section{Absorption and completion of the proof}\label{sec:completion}

We apply the stopping-tree estimate simultaneously in the three
grids $\D^0,\D^1,\D^2$. An interval $P\in\Pfamily_h^{\D^j}$ may
be assigned to an interval $Q(P)$ in another grid. Summing the three main quantities
allows such transfers to be absorbed in one inequality.

The families $\Pfamily_h^{\D^j}$ below are those defined in
\eqref{eq:boundary-family}. The threshold $T$ will be chosen
only after the coefficient of the absorbable term is known.

\begin{proposition}\label{prop:boundary-estimate}
Let $w$ satisfy the $C_p$ condition with exponent
$r>1$ and constant $K$.
For every $T\geq2$ and every minimum length $h=2^k$,
\[
 \begin{split}
 \sum_{j=0}^2\sum_{P\in\Pfamily_h^{\D^j}}
          \omega_P^{p-1}f_Pw(P)
 &\leq C_{p,r,K,T}\int_{\R}(Mf)^pw\\
 &\quad+\frac{C_pK}{T}
                 \sum_{j=0}^2\mathcal A_h^{\D^j}.
 \end{split}
\]
In particular, the constant $C_p$ in the second term is independent
of $T$, $r$, and $K$.
\end{proposition}

\begin{proof}
Recall that $\Pfamily_h^{\D^j}$ consists of intervals $P\in\D^j$
with $|P|\geq2h$ such that $J_P\leq cw_P$, while at least one
of the two dyadic children $I$ of $P$ satisfies $J_I>cw_I$.
For each such $P$, \eqref{eq:weight-drop} gives
\[
 \min(w_{P_-},w_{P_+})\leq\tau w_P.
\]
Thus $P$ belongs to $\Ofamily_\tau^{\D^j}$, the family of intervals
with a fixed relative drop in the average of $w$ on one child.

We first estimate the contribution of intervals with $R_P\leq T$.
Using $f_P\leq\omega_P$ and Corollary~\ref{cor:bounded-R}, we obtain
\[
 \sum_{j=0}^2
 \sum_{\substack{P\in\Pfamily_h^{\D^j}\\R_P\leq T}}
          \omega_P^{p-1}f_Pw(P)
 \leq C_{p,r}\log T\int_{\R}(Mf)^pw.
\]
It remains to estimate the intervals $P\in\Pfamily_h^{\D^j}$ with $R_P>T$.
For each such $P$, let $Q=Q(P)$ be the assigned interval, and
interpret $J_Q$ and $\widehat Q$ in the grid $\D^{j(P)}$ containing
$Q$. There are three possibilities.

\medskip
\noindent\emph{Small tail at $Q$: $J_Q\leq cw_Q$.}
The interval $Q$ is faithful with constant $F_p$, so the $C_p$ condition
gives
\[
 R_Qw_Q
 \leq\frac K{|Q|}\int_{\R}(M\one_Q)^pw
 \leq KF_pJ_Q
 \leq KF_pcw_Q.
\]
Put $T_0=\max(2,KF_pc)$. For each $k\in\{0,1,2\}$, apply
Lemma~\ref{lem:entropy-transfer} to the intervals $P$ under
consideration for which $j(P)=k$, with $\D_*=\D^k$,
$N=3$, $a=8$, and $T=T_0$.
Its hypotheses hold because $P\subset Q$, $|Q|=8|P|$,
and $R_Q\leq T_0$. Summing over $k=0,1,2$ and
using $f_P\leq\omega_P$, we obtain a total contribution bounded by
\[
 C_{p,r,K}\int_{\R}(Mf)^pw.
\]
The ratio $R_P$ may be arbitrarily large here: the entropy
estimate is applied on $Q(P)$, where the ratio is bounded.

\medskip
\noindent\emph{Large tail and rapid growth at $Q$:
$J_Q>cw_Q$ and $\omega_Q\geq\frac32\omega_{\widehat Q}$.}
For each $j\in\{0,1,2\}$, Lemma~\ref{lem:transferred-growth}
controls the sum of $\omega_P^pw(P)$ over these intervals
$P\in\Pfamily_h^{\D^j}$.
Again $f_P\leq\omega_P$, so the total contribution of this case
is at most
\[
 C_p\int_{\R}(Mf)^pw.
\]

\medskip
\noindent\emph{Large tail and slow growth at $Q$:
$J_Q>cw_Q$ and $\omega_Q<\frac32\omega_{\widehat Q}$.}
The estimate \eqref{eq:selection-R} gives
$w_P\leq C_pKJ_Q/R_P$. Since $P\subset Q$ and $f\geq0$,
we also have $\int_Pf\leq\int_Qf$. Using these two estimates
and \eqref{eq:selection-f}, we obtain
\begin{equation}\label{eq:small-coefficient}
 \begin{split}
 \omega_P^{p-1}f_Pw(P)
 &=\omega_P^{p-1}\left(\int_Pf\right)w_P\\
 &\leq\frac{8^{p-1}C_pK}{R_P}
             \omega_Q^{p-1}\left(\int_Qf\right)J_Q\\
 &=\frac{8^{p-1}C_pK}{R_P}
             \omega_Q^{p-1}f_QJ_Q|Q|\\
 &\leq\frac{8^{p-1}C_pK}{T}
             \omega_Q^{p-1}f_QJ_Q|Q|.
 \end{split}
\end{equation}
The term indexed by $Q$ belongs to $\mathcal A_h^{\D^{j(P)}}$.
Indeed, the tail and growth conditions are precisely those in
\eqref{eq:A}, and
$|Q|=8|P|\geq16h$ preserves the lower cutoff.

For each fixed $Q$ and each grid $\D^j$, at most eight intervals
$P\in\D^j$ can satisfy $Q(P)=Q$, because they all have length
$|Q|/8$ and lie inside $Q$. The total multiplicity of $Q$ is
thus at most $24$. Summing \eqref{eq:small-coefficient} and
absorbing the factor $24\cdot8^{p-1}$ into $C_p$ bounds this case by
\[
 \frac{C_pK}{T}\sum_{j=0}^2\mathcal A_h^{\D^j}.
\]
Combining the low-ratio estimate with the three cases proves
the proposition. All counting estimates can first
be applied to finite subfamilies of the original intervals and then extended by
monotone convergence.
\end{proof}

The factor $f_P$ in estimate \eqref{eq:entrance-cost} is used in
\eqref{eq:small-coefficient} through the inequality
$f_P|P|\leq f_Q|Q|$. Replacing the boundary cost by
$\omega_P^pw(P)$ would lose this comparison with the summand
in $\mathcal A_h$. This is why we retained the ordinary average
throughout the stopping-tree argument.

\begin{proof}[Proof of Theorem~\ref{thm:main}]
Sum \eqref{eq:tree} over the three fixed grids and apply
Proposition~\ref{prop:boundary-estimate}. For every $T\geq2$,
\[
 \sum_{j=0}^2\mathcal A_h^{\D^j}
 \leq C_{p,r,K,T}\int_{\R}(Mf)^pw
       +\frac{C'_pK}{T}\sum_{j=0}^2\mathcal A_h^{\D^j}.
\]
The sums on both sides are finite, as verified after
\eqref{eq:A}. Choose $T\geq\max(2,2C'_pK)$ and absorb the last
term. We conclude that
\begin{equation}\label{eq:three-sums}
 \sum_{j=0}^2\mathcal A_h^{\D^j}
 \lesssim_{p,r,K}\int_{\R}(Mf)^pw,
\end{equation}
uniformly in the minimum length $h$.

Let $\Sfam$ be the original finite sparse family. Apply
Proposition~\ref{prop:sparse-transfer} to obtain
$\eta/8$-sparse families $\Sfam_j\subset\D^j$ of faithful
intervals. Choose $h$ no larger than the smallest interval in
these families. By \eqref{eq:sparse-transfer},
Proposition~\ref{prop:sparse-reduction}, and \eqref{eq:three-sums},
\begin{align*}
 \int_{\R}(A_{\Sfam}f)^pw
 &\lesssim_p
       \sum_{j=0}^2\int_{\R}(A_{\Sfam_j}f)^pw\\
 &\lesssim_{p,r,K,\eta}
       \int_{\R}(Mf)^pw+\sum_{j=0}^2\mathcal A_h^{\D^j}\\
 &\lesssim_{p,r,K,\eta}\int_{\R}(Mf)^pw.
\end{align*}
The constants are independent of the finite family. The
preliminary regularization of $w$ and monotone convergence in
$f$ and in the sparse family now prove Theorem~\ref{thm:main}
in full generality.
\end{proof}

\section{Concluding remarks}\label{sec:conclusion}

\begin{remark}[The Fefferman--Stein inequality]\label{rem:FS}
Yabuta~\cite{Yabuta} proved that $C_p$ is necessary for the weighted
Fefferman--Stein inequality and that $C_{p+\varepsilon}$, for some
$\varepsilon>0$, is sufficient. Together with Yabuta's necessity
result, Theorem~\ref{thm:main} shows that the $C_p$ condition is
necessary and sufficient for the weighted Fefferman--Stein
inequality.

\Needspace{8\baselineskip}
\begin{corollary}\label{cor:FS}
Let $1<p<\infty$ and $w\in C_p$. If $f\in L^1_{\mathrm{loc}}(\R)$
and $|\{x:|f(x)|>t\}|<\infty$ for every $t>0$, then
\[
 \int_{\R}|f|^p w\leq C\int_{\R}(f^\#)^p w,
\]
where $f^\#$ is the Fefferman--Stein sharp maximal function.
\end{corollary}

\begin{proof}
For a fixed sufficiently small $\lambda>0$, the sparse
oscillation estimate and the local sharp maximal function
$M_\lambda^\# f$ satisfy
\[
 |f|\lesssim A_{\Sfam}(M_\lambda^\# f),
 \qquad M(M_\lambda^\# f)\lesssim f^\#;
\]
see \cite[Section~2.2 and Theorem~2.3]{LernerWeak} for the definition of $M_\lambda^\# f$ and the above inequalities.
Apply Theorem~\ref{thm:main}, first to bounded compactly supported
truncations of $M_\lambda^\# f$ and then pass to the limit.
This gives
\[
 \int_{\R}|f|^p w
 \lesssim\int_{\R}\bigl[M(M_\lambda^\# f)\bigr]^p w
 \lesssim\int_{\R}(f^\#)^p w.\qedhere
\]
\end{proof}
\end{remark}

\begin{remark}[The $SC_p$ condition]\label{rem:SCp}
There is also a consequence for the class $SC_p$ introduced in
\cite[Definition~1.1]{LernerWeak}. In one dimension, a weight $w$
belongs to $SC_p$ if there exist $r>1$ and $C<\infty$ such that,
for every finite family of pairwise disjoint bounded intervals
$I_1,\ldots,I_N$,
\[
 \sum_{j=1}^N
       \left(\avg{w^r}{I_j}\right)^{1/r}|I_j|
 \leq C\int_{\R}
       \bigl(M\one_{\bigcup_{j=1}^N I_j}\bigr)^p w.
\]
The constants are independent of the family and its cardinality.
The formulation for countable families follows by monotone
convergence.

Taking a family consisting of a single interval and dividing by
its length gives exactly the $C_p$ condition.
Thus $SC_p\subset C_p$ follows directly from the definitions.
The additional requirement in $SC_p$ is to control the sum of
the local $L^r$ averages for an arbitrary disjoint family by
one tail, formed from the indicator of its union.

The reverse inclusion follows from the Fefferman--Stein
inequality proved above. If $w\in C_p$, then for every $t>0$
and every $f$ as in Corollary~\ref{cor:FS}, Chebyshev's inequality
gives
\[
 t^p w(\{x:|f(x)|>t\})
 \leq\int_{\R}|f|^pw
 \leq C\int_{\R}(f^\#)^pw.
\]
This is the weighted weak type Fefferman--Stein inequality.
By \cite[Theorem~1.3]{LernerWeak}, that inequality is equivalent
to $w\in SC_p$. Hence $C_p\subset SC_p$, and therefore
\[
 SC_p=C_p\qquad(1<p<\infty).
\]
In particular, the apparently stronger requirement for disjoint
families imposes no additional restriction on the weight.
\end{remark}

\begin{remark}[Higher dimensions]\label{rem:higher-dimensions}
Theorem~\ref{thm:main}, Corollary~\ref{cor:FS}, and the identity
$SC_p=C_p$ extend to $\R^n$, with intervals replaced by cubes
whose sides are parallel to the coordinate axes.
We indicate the adjustments to the proof.

A dyadic cube has $2^n$ children. The appropriate tails are
\[
 q=2^{-n(p-1)},\qquad
 J_Q^{\D}=\sum_{k\geq0}q^k w_{Q^{(k)}}.
\]
They satisfy the same recursion $J_R=qJ_Q+w_R$ for every child
$R$ of $Q$. In the stopping-tree calculation, the factor $2$
coming from the number of children becomes $2^n$; we choose
$c$ so that $q+2^n/c<(2/3)^{p-1}$. The resulting positive
coefficient and the weight drop are therefore preserved.
The entropy argument uses convexity of $t\log t$ for the
$2^n$ child averages: a fixed relative drop on one child still
gives an increment bounded below by a positive multiple of
$w(Q)$, and the increments telescope as before. The maximal-average
comparisons and the embedding argument also extend, with
constants depending on $n$.

For the geometric step, take the $3^n$ product grids formed
from the three one-dimensional grids in Section~\ref{sec:grids}.
The two-scale argument in Lemma~\ref{lem:selection} applies
coordinate by coordinate. In each coordinate, at least two of
the three grids are suitable at each of the two scales, so a
common choice exists. The resulting product grid gives an
assigned cube $Q(P)\supset P$ with
$\ell(Q(P))=8\ell(P)$, where $\ell$ denotes side length,
and the same tail comparisons, with constants depending on
$n$ and $p$. At most $8^n$ cubes from any one original grid
can be assigned to the same cube.

In Lemma~\ref{lem:transferred-growth}, the sixth ancestor is
replaced by the $m$th ancestor, where $m=m(n)$ is large enough
that $(1+2^{4-m})^n<3/2$. Indeed, when
$\omega_{Q(P)}\geq\frac32\omega_{\widehat{Q(P)}}$, the same
enlargement of a cube containing $P^{(m)}$ gives
\[
 \omega_{P^{(m)}}
 \leq (1+2^{4-m})^n\omega_{\widehat{Q(P)}}
 \leq\tfrac23(1+2^{4-m})^n\omega_P.
\]
The coefficient of $\omega_P$ is strictly less than one, so
Lemma~\ref{lem:growth} again applies. These changes preserve the
absorption argument and prove the sparse estimate in every
dimension. The proofs in Remarks~\ref{rem:FS} and~\ref{rem:SCp}
then apply to cubes as well.

Finally, standard sparse domination, see
\cite[Theorem~2.2]{LernerWeak}, shows that $C_p$ is sufficient
for the corresponding inequality for Calder\'on--Zygmund
operators. More precisely, if $1<p<\infty$, $w\in C_p(\R^n)$,
and $T$ is an $L^2$-bounded Calder\'on--Zygmund operator with
a kernel satisfying the standard size and Dini continuity
conditions, then
\[
 \int_{\R^n}|Tf|^p w
 \leq C\int_{\R^n}(Mf)^p w
\]
for every bounded compactly supported $f$. The constant depends
only on $n,p$, the parameters of the $C_p$ condition, and the
operator bounds and kernel constants of $T$.
\end{remark}

\appendix
\section{Why the sparse sum cannot be estimated directly}\label{app:sparse-sum}

The $C_p$ condition does not, in general, allow us to bound the sum
on the right-hand side of \eqref{eq:sparse-r} by a constant multiple
of $\int_{\R}(Mf)^p w$. The following example shows that this sum
may be infinite even when the maximal-function integral is finite.
Moreover, choosing the exponent $r$ closer to one does not remove
this obstruction.

\begin{proposition}\label{prop:sparse-counterexample}
Let $1<p<\infty$ and define
\[
 w=\sum_{k=2}^{\infty}
 \frac{2^{k(p-1)}}{k^{1-1/p}(\log k)^{2-1/p}}
 \one_{[2^k,\,2^k+2^k/(k\log k)^{1/p}]}.
\]
Then $w\in C_p$, and the defining $C_p$ inequality holds for every
$r>1$. There are a dyadic grid $\D$ and a $3/4$-sparse family
$\Sfam\subset\D$ such that, for $f=\one_{[0,1]}$,
\[
 \int_{\R}(Mf)^p w<\infty,
 \qquad
 \sum_{I\in\Sfam}f_I^p
       \left(\avg{w^r}{I}\right)^{1/r}|I|=\infty
 \quad\text{for every }r>1.
\]
In particular, for this fixed weight and any $r>1$, the right-hand
side of \eqref{eq:sparse-r} cannot be bounded by a constant multiple
of $\int_{\R}(Mf)^p w$ uniformly over finite sparse families.
\end{proposition}

\begin{proof}
Write
\[
 \ell_k=\frac{2^k}{(k\log k)^{1/p}},\qquad
 E_k=[2^k,2^k+\ell_k],\qquad
 a_k=\frac{2^{k(p-1)}}{k^{1-1/p}(\log k)^{2-1/p}},
\]
so that $w=\sum_{k\geq2}a_k\one_{E_k}$. We have $\ell_k<2^k$,
so these intervals are disjoint. The two identities
\[
 a_k=\frac{\ell_k^{p-1}}{\log k},
 \qquad
 \frac{a_k\ell_k}{2^{kp}}=\frac{1}{k(\log k)^2}
\]
explain the choice of their lengths and heights. The weight is
locally bounded, and the second identity gives
\[
 \int_4^{\infty}\frac{w(x)}{x^p}\,dx
 \leq\sum_{k=2}^{\infty}\frac{1}{k(\log k)^2}<\infty.
\]
Consequently, $\int_{\R}(M\one_I)^p w$ is finite for every bounded
interval $I$: on a bounded set this follows from local integrability,
and for large positive $x$ one has $M\one_I(x)\lesssim_I x^{-1}$.

\medskip
\noindent\emph{Verification of the $C_p$ condition.}
We prove the stronger estimate
\begin{equation}\label{eq:appendix-supremum}
 \operatorname*{ess\,sup}_{I}w
 \lesssim_p \frac{1}{|I|}\int_{\R}(M\one_I)^p w
 \qquad\text{for every bounded interval }I.
\end{equation}
Fix $I$, put $L=|I|$, and consider any index $n\geq2$ for which
$|I\cap E_n|>0$. It is enough to bound $a_n$ by the right-hand side
of \eqref{eq:appendix-supremum}, with a constant independent of $I$
and $n$.

First suppose $L\leq\ell_n$. Choose $z\in I\cap E_n$. Since
$E_n$ has length at least $L$, there is a subinterval $F\subset E_n$
of length $L$ such that $|x-z|\leq L$ for $x\in F$. An interval
containing both $I$ and such an $x$ has length at most $2L$.
Therefore $M\one_I(x)\geq1/2$ on $F$, and
\[
 \frac1L\int_{\R}(M\one_I)^p w
 \geq\frac{2^{-p}}L\int_Fw=2^{-p}a_n.
\]

Now suppose $L>\ell_n$. Let $m\geq n$ be the smallest integer
such that $L\leq2^m$. Because $I$ meets $E_n\subset[2^n,2^{n+1}]$,
we have
\[
 I\subset[-2^m,3\cdot2^m].
\]
If $k\geq m+3$ and $x\in E_k$, an interval containing $I$ and
$x$ has length at most $3\cdot2^k$. Hence
\[
 M\one_I(x)\geq\frac{L}{3\cdot2^k},
\]
and summing the contributions of these distant intervals gives
\begin{align*}
 \frac1L\int_{\R}(M\one_I)^p w
 &\geq 3^{-p}L^{p-1}
        \sum_{k=m+3}^{\infty}\frac{a_k\ell_k}{2^{kp}}\\
 &=3^{-p}L^{p-1}
        \sum_{k=m+3}^{\infty}\frac{1}{k(\log k)^2}
 \gtrsim_p\frac{L^{p-1}}{\log m}.
\end{align*}
The last inequality follows by comparison with
$\int_{m+3}^{\infty}dt/(t(\log t)^2)=1/\log(m+3)$.

If $m=n$, this lower bound is at least a constant multiple of
$\ell_n^{p-1}/\log n=a_n$. If $m>n$, the choice of $m$ gives
$L>2^{m-1}$. We use the elementary comparison
\[
 \frac{2^{n(p-1)}}{\log n}
 \lesssim_p \frac{2^{m(p-1)}}{\log m}
 \qquad(m\geq n\geq2).
\]
Indeed, writing $d=m-n$, the ratio of the left-hand side to the
right-hand side is at most
$(1+d)2^{-d(p-1)}$, which is bounded in terms of $p$; here
$\log m/\log n\leq1+d$ follows from
$\log(m/n)\leq d/n$ and $n\log n\geq1$.
It follows that
\[
 \frac1L\int_{\R}(M\one_I)^p w
 \gtrsim_p\frac{2^{m(p-1)}}{\log m}
 \gtrsim_p\frac{2^{n(p-1)}}{\log n}
 \geq a_n.
\]
Thus in every case the tail controls $a_n$. Taking the maximum over
the finitely many intervals $E_n$ meeting $I$ in positive measure
proves \eqref{eq:appendix-supremum}; if there are none, its left-hand
side is zero. Finally,
\[
 \left(\avg{w^r}{I}\right)^{1/r}
 \leq\operatorname*{ess\,sup}_{I}w,
\]
so \eqref{eq:appendix-supremum} proves the $C_p$ condition for every
$r>1$, with a constant depending only on $p$.

\medskip
\noindent\emph{Failure of the bound for the sparse sum.}
Set $f=\one_{[0,1]}$. Since $Mf(x)=1/x$ for $x>1$ and $w$ is
supported in $[4,\infty)$, the preceding calculation gives
\[
 \int_{\R}(Mf)^p w
 =\sum_{k=2}^{\infty}a_k\int_{E_k}x^{-p}\,dx
 \leq\sum_{k=2}^{\infty}\frac{1}{k(\log k)^2}<\infty.
\]

For $j\geq1$, let
\[
 I_j=[-4^j,2\cdot4^j),\qquad
 \Sfam=\{I_j:j\geq1\}.
\]
Each $I_j$ is the second of the four dyadic grandchildren of
$I_{j+1}$. The collection of all dyadic subintervals of the $I_j$
therefore forms a dyadic grid. Since $I_j\uparrow\R$, this grid
satisfies the definition in Section~\ref{sec:prelim}. The sets
$I_1$ and $I_j\setminus I_{j-1}$, $j\geq2$, witness that $\Sfam$
is $3/4$-sparse.

We have $E_{2j}\subset I_j$, $|I_j|=3\cdot4^j$, and
$f_{I_j}=|I_j|^{-1}$. Fix $r>1$ and put
\[
 \beta=\frac{p-1+1/r}{p}=1-\frac{1-1/r}{p}<1.
\]
Keeping only the contribution of $E_{2j}$ in the $L^r$ average,
we obtain
\begin{align*}
 f_{I_j}^p\left(\avg{w^r}{I_j}\right)^{1/r}|I_j|
 &\geq |I_j|^{1-p-1/r}a_{2j}\ell_{2j}^{1/r}\\
 &=\frac{3^{1-p-1/r}}
        {(2j)^\beta(\log(2j))^{1+\beta}}.
\end{align*}
The series of these lower bounds diverges because $\beta<1$.
Thus the sparse sum is infinite for every $r>1$. Taking the finite
subfamilies $\{I_1,\ldots,I_N\}$ and letting $N\to\infty$ proves
the asserted failure of a uniform estimate for finite families.
\end{proof}

\medskip
\noindent\emph{The role of the rapid-growth subfamily.}
Every interval in this example belongs to the rapid-growth subfamily.
Both $I_j$ and its dyadic parent contain the support of $f$, so
\[
 \omega_{I_j}=\frac1{|I_j|},\qquad
 \omega_{\widehat I_j}=\frac1{2|I_j|}.
\]
In particular, $\omega_{I_j}=2\omega_{\widehat I_j}$, and none of
these intervals satisfies the slow-growth restriction in
\eqref{eq:A}. Their contribution to the operator itself is harmless.
Indeed, if $j$ is the first index for which $x\in I_j$, then
\[
 A_{\Sfam}f(x)=\sum_{k=j}^{\infty}\frac1{3\cdot4^k}
 =\frac43\frac1{|I_j|}\leq\frac43Mf(x).
\]
Thus the operator is pointwise controlled by $Mf$, whereas the sum
in \eqref{eq:sparse-r} diverges. Proposition~\ref{prop:sparse-reduction}
therefore removes this subfamily before applying
Lemma~\ref{lem:sparse} to the remaining part.

\vskip 3mm
\noindent

{\bf Acknowledgements and use of AI.}
AL was supported by ISF grant no. 807/26. SO was partially supported by CONICET, grant  KE2.11220250100060CO.

The first version of this paper was written jointly with Fedor Nazarov, who was originally expected to be a coauthor but ultimately chose not to be.
We thank Fedor Nazarov for numerous discussions, new ideas, and for his interest in this work. That version was prepared without the use of AI and included a complete proof of the main results. 
The proof already used Lemma~3.1, after first treating separately the intervals on which the maximal averages of \(f\) grow rapidly. The estimate for the remaining sum relied on Bellman-type and probabilistic arguments.

For the present version, ChatGPT was used to develop an alternative treatment of the remaining sum, replacing the Bellman-type and probabilistic arguments with the deterministic arguments presented in Sections 4–7. The authors have checked the resulting proof and take full responsibility for the mathematical content of this paper.

\end{document}